\documentclass[10pt,a4paper]{article}
\usepackage[utf8]{inputenc}
\usepackage[english]{babel}
\usepackage{amsmath}
\usepackage{amsfonts}
\usepackage{amssymb}
\usepackage{graphicx}
\usepackage{appendix}
\usepackage{float}

\begin{document}

\title{Natural Factor based Solvers}
\author{
O. Andrés Cuervo\thanks{School of Engineering, Science and Technology, Universidad del Rosario, Bogot\'a, Colombia. {\tt omar.cuervo@urosario.edu.co}}, Juan Galvis\thanks{Departamento de Matem\'aticas, Universidad Nacional de Colombia, Bogot\'a, Colombia.{\tt jcgalvisa@unal.edu.co}}, Marcus Sarkis\thanks{Department of Mathematical Sciences, Worcester Polytechnic 
Institute Worcester USA} 
}

\maketitle
\begin{abstract}
We consider parametric families of partial differential 
equations--PDEs where the parameter $\kappa$ modifies  only the (1,1) block of a saddle point matrix product of a discretization below. The main goal is
to develop an algorithm that removes, as much as possible, the dependence
of iterative solvers on the parameter $\kappa$. The algorithm we propose
requires only one matrix factorization which does not depend on $\kappa$,
therefore, allows to reuse it for solving very fast a large number of
discrete PDEs for different $\kappa$ and forcing terms. The design of
the proposed algorithm is motivated by previous works on natural factor
of formulation of the stiffness matrices and their stable numerical solvers.
As an application, in two dimensions, we consider an iterative preconditioned
solver  based on the null space of Crouzeix-Raviart discrete gradient
represented as the discrete curl of $P_1$ conforming finite element functions. 
For the numerical examples, we consider the case of random coefficient pressure equation where the permeability is modeled by an stochastic process. We note
that contrarily from recycling Krylov subspace techniques, the proposed
algorithm does not require fixed forcing terms.

\end{abstract}
\section{Introduction}

The general form of a saddle point system of linear equations we consider is
\begin{equation}\label{saddle_point_problem}
  \left[ \begin{array}{cc}
D(\kappa)^{-1} & G \\
G^T & 0 
\end{array} \right] \left[ \begin{array}{c}
     q  \\
      u
\end{array} \right]=\left[  \begin{array}{c}
     r  \\
      b
\end{array} \right],
\end{equation} 
where the matrix $D$ is symmetric positive definite. This form is standard
in the formulation of mixed finite elements. What is not very
well-known, as pointed out by Argyris and Br{\o}nlund \cite{argyris1975natural},
is that classical conforming and nonconforming finite element methods -- FEMs
can also be written in the  form (\ref{saddle_point_problem}) with $r=0$;
see Section \ref{CRdisc} for the case of
Crouzeix-Raviart FEM and Section \label{AL} for $P_1$ conforming FEM. We show
that the stiffness matrix,  associated
to the Crouzeix-Raviart FEM element discretization for
the PDE (\ref{eq:hetdiff}) with isotropic coefficients
$\kappa(x)$, has the the natural factor  of the form
$A_{CR} = G_{CR}^T D(\kappa) G_{CR}$, where $G_{CR}$ is the discrete gradient
(not affected by the parameter $\kappa$) and  $D(\kappa)$ is a
diagonal matrix with entries depending of the integration
of $\kappa$ in each element, hence, it is easy to update the
natural factor if $\kappa$ is modified. Due to the 
superior numerical stability with respect to roundoff errors when 
operating with $G^T$, $D(\kappa)$ and $G$ rather than the assembled
stiffness matrix, several works  \cite{vavasis1996stable, vavasis1994stable, rees2018comparative, drmavc2021numerical} were dedicated in solving the saddle point problem 
(\ref{saddle_point_problem}) or associated SVD and diagonalization. In Sections \ref{NullSpace} and \ref{RangeSpace}  we review some aspects of these works. The methods start by representing
$q$ on the range of the matrix $[G\,\, Z]$ where $Z$ is such that
$Q = [G \,\, Z]$ is a square invertible matrix; two common choices of $Z$
are $Z^T G = 0$ or $Z^T D^{-1} G = 0$. These works generate very stable
algorithms for ill-conditioned $\kappa$, however,  they do not remove the
dependence on $D$ of the factorizations, hence, they do not fit our
goal of reusing the same factorization for different values of $\kappa$.
In Section \ref{auxpro} we propose our method, we first use discrete
Hodge Laplacian ideas to choose $Z=\tilde{C}_L$ as the curl of $P_1$ conforming
piecewise linear basis functions, hence $G^T_{CR} \tilde{C}_L = 0$. Then we consider the coupled system 
\[
(\text{grad}\, u_{CR} + \text{curl}\, w_{P_1}, 
\kappa\, (\text{grad}\, v_{CR} + \text{curl}\, v_{P_1}))_{L^2(\Omega)} 
\]
as a preconditioner for the uncoupled system 
\[
(\text{grad}\, u_{CR}, \kappa\, \text{grad}\, v_{CR})_\kappa + (\text{curl}\, w_{P_1},\kappa\, \text{curl}\  v_{P_1})_{L^2(\Omega)}.
  \]

\section{Crouzeix-Raviart nonconforming finite elements} \label{CRdisc}
Consider the heterogeneous diffusion equation

\begin{equation}\label{eq:hetdiff}
 \left\lbrace \begin{array}{rclr}
-\partial_1(\kappa(x)\partial_1 u(x))-\partial_2(\kappa(x)\partial_2 u(x)) & = & f(x),  & \ \ \ x\in \Omega, \\
u(x) & = & 0, & \ \ \ x \in \partial \Omega,
\end{array} \right. 
\end{equation}
where $\Omega \subseteq \mathbb{R}^2$ and  $\kappa: \Omega \rightarrow \mathbb{R}^+$, $f: \Omega \rightarrow \mathbb{R}$ are given. 

In particular, in the target application $\kappa(x)$ is a random field that describes the permeability and allows modeling the lack of data and uncertainties of the problem (e.g., subsurface flow). The forcing term $f$ may also be a random field. In general, in many practical situations we must solve \eqref{eq:hetdiff} for a large family of coefficients $\kappa$ and forcing terms $f$. See Section \ref{application}.

Let us introduce a triangulation $\mathcal{T}^h$ of $\Omega$. Discretize \eqref{eq:hetdiff} by the Crouzeix-Raviart (CR) non-conforming finite element space.   Define the CR space $\widetilde{V}^{CR}$ as the space of all piecewise linear functions with respect to $\mathcal{T}^h$ that are continuous at interior edges midpoints. The degrees of freedom are located in the midpoint of the edges of
$\mathcal{T}^h$. Let $V_{CR}\subseteq \widetilde{V}^{CR}$ the subspace of functions in $\widetilde{V}^{CR}$ with zero value at the midpoint of boundary edges. The approximation $u_{CR}
\in V_{CR}$ of the solution of (\ref{eq:hetdiff}) is the solution of 
\begin{equation*}
\sum_{T \in \mathcal{T}^h} \int_T \kappa(x)(\partial_1u_{CR} (x)\partial_1v (x) + \partial_2u_{CR}(x)\partial_2v(x)) dx = \int_\Omega f(x)v(x) dx,
\end{equation*}
for all $v \in V^{CR}.$ The  linear system of the $CR$ approximation is given by 
\begin{equation}\label{eq:ACR}
A_{CR} u_{CR} = b_{CR},
\end{equation} where 
$
A_{CR}=\left[ a_{ei}^{CR} \right]_{e,i=1}^{N_e}$ and $ b_{CR}=\left[ b_e \right]_{e=1}^{N_e}.
$
Here, $N_e$ denotes the number of interior edges of $\mathcal{T}^h$, 
$
b_e =\int_\Omega f(x) \phi_e^{CR}(x) dx
$ and \\
$$
a_{ei}^{CR}=\sum_{T \in \mathcal{T}^h} \int_T \kappa(x)\Big(\partial_1\phi^{CR}_e (x)\partial_1\phi_i^{CR} (x) +\partial_2\phi^{CR}_e (x)\partial_2\phi_i^{CR} (x)\Big) dx.
$$

Let $x_T$ denote the barycenter of triangle  $T \in  \mathcal{T}^h$. Piecewise gradients of functions in $V^{CR}$ are piecewise constant vector functions and then 
\begin{equation}
a_{ei}^{CR}=\sum_{T \in \mathcal{T}^h} \kappa_T \vert T \vert \partial_1 \phi_e(x_T) \partial_1 \phi_i(x_T) + \sum_{T \in \mathcal{T}^h} \kappa_T \vert T \vert \partial_2 \phi_e(x_T) \partial_2 \phi_i(x_T)
\end{equation}
where $\kappa_T$ is the average value of $\kappa(x)$ in $T$.  Therefore, we can write (see \cite{argyris1975natural}) 
\begin{equation*}
A_{CR} = G^T_{CR} D G_{CR} = G_{CR,1}^T D_1 G_{CR,1} +G_{CR_2}^T D_2 G_{CR,2}
\end{equation*} where $G_{CR,l}= \left[ g_{e,T}^{CR,l} \right]_{N_T \times N_e} = \left[ \sqrt{\vert T \vert} \partial_l \phi^{CR}_e(x_T) \right]_{N_T \times N_e},$  
and $N_T$ denotes the number of triangles in $\mathcal{T}^h$  and  $l=1,2$. Furthermore, write,
\begin{equation}\label{def:Dkappa}
D_l= \mbox{diag}(\kappa_T)_{T \in \mathcal{T}^h},  \ \ D=\mbox{diag}(D_1, D_2)  \ \ \mbox{and} \ \ G_{CR}=\left[ \begin{array}{c}
G_{CR,1} \\
G_{CR,2}
\end{array} \right]_{2N_T\times N_e}.
\end{equation} 
We can write the matrix formulalation as 
\begin{equation}\label{eq:Schur1}
    G_{CR}^TDG_{CR}u_{CR}=b_{CR}.
\end{equation}
We see that problem (\ref{eq:Schur1}) is the Schur complement of the saddle point problem
\begin{equation}\label{Schur1}
 \left[ \begin{array}{cc}
D^{-1} & G_{CR} \\
G_{CR}^T & 0 
\end{array} \right] \left[ \begin{array}{c}
     q  \\
      u_{CR}
\end{array} \right]=\left[  \begin{array}{c}
     0  \\
     - b_{CR}
\end{array} \right].
\end{equation}

\section{Conforming finite elements $P_1$} \label{AL}

Let 
$
\widetilde{V}^L = P_1(\mathcal{T}^h)=\{ v:\Omega\rightarrow \mathbb{R} \vert \ \  v\vert_{T} \mbox{ is linear for all } T \in \mathcal{T}^h  \} \cap C^0(D). 
$
The space $\widetilde{V}^L$ has a base $ \{ \varphi_i^L \}_{i=1}^{\widetilde{N}_{v}} $, where $\widetilde{N}_v$ is the number of vertices and $\varphi_i^L$ is the function that takes value 1 at the $i - th$ node and 0 at the other nodes. 
Also define $V^L=\widetilde{V}^L \cap H_0^1(\Omega)$ and $N_{v}$ the number of interior vertices. \\
The approximation $u_L$ of the solution of (\ref{eq:hetdiff}) is: find $u_L \in V^L$ such that  
\begin{equation*}
 \int_\Omega \kappa(x)(\partial_1u_{L} (x)\partial_1v (x) + \partial_2u_{L}(x)\partial_2v(x)) dx = \int_\Omega  f(x)v(x) dx,
\end{equation*} for all $v \in V^L$, with matrix form $$A_Lu_L=b_L,$$ where, $A_L=[a_{ij}^L]_{i,j=1}^{N_{v}}$ and $b_L=[b_i^L]_{i=1}^{N_v}$ with 
$
b_i^L =\int_D f(x) \varphi_i^{L}(x) dx
$ and 
\begin{equation*}
a_{ij}^{L}= \int_\Omega \kappa(x)\Big(\partial_1\varphi^{L}_i (x)\partial_1\varphi_j^{L} (x) +\partial_2\varphi^{L}_i (x)\partial_2\varphi_j^{L} (x)\Big) dx.
\end{equation*}
As before, we have (see \cite{argyris1975natural}) 
\begin{equation*}
A_L=G_L^T D G_L=G_{L,1}^TD_1G_{L,1}+G_{L,2}^TD_2G_{L,2}  
\end{equation*} where 
 $G_{L,l}= \left[ g_{e,v}^{L,l} \right]_{N_T \times N_v} =\left[ \sqrt{\vert T \vert } \partial_l \varphi_v^L(x_T)  \right]_{N_T \times N_v}$ and $G_L=\left[ \begin{array}{c}
G_{L,1} \\
G_{L,2}
\end{array} \right]_{2N_T\times N_v}.$
We can write the matrix formulation as 
\begin{equation}\label{Schur2}
    G_L^TDG_Lu_L=b_L, 
\end{equation}
and the corresponding saddle point problem is 
\begin{equation*}
 \left[ \begin{array}{cc}
D^{-1} & G_{L} \\
G_{L}^T & 0 
\end{array} \right] \left[ \begin{array}{c}
     q  \\
      u_{L}
\end{array} \right]=\left[  \begin{array}{c}
     0  \\
      -b_{L}
\end{array} \right].
\end{equation*}

\section{The null space method} \label{NullSpace} 
A method for solving the saddle point problem (\ref{saddle_point_problem}) is called the null space method, see \cite{rees2018comparative}. We split (\ref{saddle_point_problem}) into two equations,
\begin{equation}\label{saddle1}
    D^{-1}q+Gu=r \quad \quad 
    G^Tq=b.
\end{equation}

The null space method consists in find $Z$  that represents the null space of $G^T$, 
$G^TZ=0$, and such that $[G \ \ Z]$ is a non-singular square matrix. Therefore, we can change variables to potentials $\chi$ and $\psi$ such that 
\begin{equation}\label{q}
    q= [G \ \ Z] \left[ \begin{array}{c}
         \chi  \\
        \psi
    \end{array} \right] = G\chi + Z \psi.
\end{equation}
From (\ref{q}) and $G^TZ=0$ we have 
$ G^Tq=G^TG\chi 
$ and from (\ref{saddle1}) we have $ b=G^TG\chi$ which gives $
  \chi=(G^TG)^{-1}b$,
that can be pre-computed. On the  other hand, from (\ref{saddle1}) and (\ref{q})  we have that 
\begin{equation*}
    D^{-1}G\chi+D^{-1}Z\psi + Gu=r   
\end{equation*}
which  gives $Z^TD^{-1}Z\psi = Z^Tr-Z^TD^{-1}G\chi$ and if we call ${c}=Z^Tr-Z^TD^{-1}G\chi$, we can write the system 
\begin{equation} \label{Zproblem} 
    Z^TD^{-1}Z\psi={c}.
\end{equation}
This is the null space system and it is similar that the Schur complement
of (\ref{saddle_point_problem}) given by $G^T D G u = -g$. 

   \section{Range null-space hybrid} \label{RangeSpace} 
To avoid solving the equation (\ref{Zproblem}), now we combine the equations (\ref{saddle1}) and (\ref{q}), we have 
$D^{-1}(G\chi + Z \psi)+Gu=r$ which gives $Z\psi+DGu=Dr-G\chi$ and it allows as to write the system  (\cite{vavasis1994stable, vavasis1996stable})
\begin{eqnarray} \label{ztilde2} 
 [DG \ \ Z] \left[ \begin{array}{c}
         u  \\
        \psi
    \end{array} \right] = Dr-G\chi.
\end{eqnarray}
We note  that the matrix $[ DG \ \ Z]$ is a square matrix and this system is  called range space scaled system. The related matrix $[G \ \ D^{-1}Z]$ is called null space scaled matrix. This algorithm is called  ``hybrid'' because
uses both the range-space and the null-space. See \cite{vavasis1994stable, vavasis1996stable}.

Alternatively, we can proceed as follows, we multiply (\ref{saddle1}) by $Z^T$
to get $Z^T(D^{-1}q+Gu)=Z^Tr$ which gives $Z^TD^{-1}q=Z^Tr$ and together with 
(\ref{q}) gives the system 
\begin{eqnarray} \label{ztilde} 
 \left[ \begin{array}{c}
         G^T  \\
        Z^TD^{-1}
    \end{array} \right] q = \left[ \begin{array}{c}
         b  \\
        Z^Tr
    \end{array} \right].
\end{eqnarray}
Note that the matrices (\ref{ztilde2}) or (\ref{ztilde})
has a dependence on $D$, however, for
numerical stability purpose is very efficient since the matrix is based on
discrete gradient times $D$ rather than the assembled second-order derivatives
with $D$. There are versions where $Z$ is replaced by $D\widetilde{Z}$,
or equivalently $G^T D \widetilde{Z}=0$, hence the matrix in (\ref{ztilde})
does not depend on $D$; unfortunately $\widetilde{Z}$ depends on $D$.

\section{An auxiliary problem and $2\times2$  systems} \label{auxpro}
Recall that for a scalar $w$,  $\overrightarrow{\mbox{curl }}w=(\partial_2 w, - \partial_1 w)$ and for a vector $\overrightarrow{q}=(q_1,q_2)$, $\mbox{curl } \overrightarrow{q}=\partial_1 q_2- \partial_2 q_1$. Consider now the elliptic equation
\[
 \left\lbrace \begin{array}{rclr}
-\mbox{curl }(\kappa(x) \overrightarrow{\mbox{curl }} w(x))& = & g(x),  & \ \ \ x\in \Omega \\
 \kappa(x) \mbox{curl } w(x) \cdot 
 \vec{ \tau} & = & 0, & \ \ \ x \in \partial \Omega
\end{array} \right.
\] 
where $\tau$ is the tangential vector on the boundary of $\Omega$. Note that we have 
$
\mbox{curl }(\kappa(x) \overrightarrow{\mbox{curl }} w(x))= -\partial_1(\kappa(x) \partial_1 w(x))-\partial_2(\kappa(x) \partial_2 w(x))
$ and
\begin{eqnarray*}
 \kappa(x) \overrightarrow{ \mbox{curl }} w(x) \cdot 
 \vec{\tau} & = & \tau_1 \kappa(x)\partial_2 w(x)- \tau_2 \kappa(x) \partial_1 w(x) \\
 & = & - n_2 \kappa(x) \partial_1 w(x) - n_1\kappa(x) \partial_2 w(x) 
 =  -\kappa (x) \nabla w(x)\cdot  \vec{\eta},
\end{eqnarray*}
 where $\vec{\eta}$ is the normal vector. We approximate this problem by conforming elements. Let 
$
\widetilde{V}^L = P_1(\mathcal{T}^h)=\{ v:\Omega\rightarrow \mathbb{R} \vert v\vert_{T} \mbox{ is linear for all } T \in \mathcal{T}^h  \} \cap C^0(\Omega).
$
The approximation of the problem above is: Find $\widetilde{w}_L\in \widetilde{V}^L$ such that 
\begin{equation*}
\int_\Omega \kappa(x) \overrightarrow{\mbox{curl }} \widetilde{w}_L(x) \cdot \overrightarrow{\mbox{curl }} v(x)  dx = \int_\Omega g(x)v(x) dx \ \ \ \ \ \mbox{ for all } v \in \widetilde{V}^L,
\end{equation*} with additional requirement that $\int_\Omega \widetilde{w}_L(x) dx =0$. The matrix form is given by 
\begin{equation*}
\widetilde{A}_L \widetilde{w}_L=\widetilde{b}_L,
\end{equation*} where $\widetilde{A}_L=[{a}^L_{ij}]_{\widetilde{N}_v \times \widetilde{N}_v}$ and $\widetilde{b}_L=[{b}_i^L]_{\widetilde{N}_v \times 1}$ with entries defined by
$
{a}_{ij}^L=\int_\Omega \kappa(x)\overrightarrow{\mbox{curl }} \varphi_i^L(x) \cdot \overrightarrow{\mbox{curl }}\varphi_j^L(x)dx$  and  ${b}_i^L=\int_\Omega g(x)\varphi_i^L(x) dx
$.
Here $\widetilde{N}_v$ is the number of vertices in $\mathcal{T}^h$. As before, we have 
\begin{equation*}
\widetilde{A}_L=\widetilde{C}_L^T D \widetilde{C}_L=\widetilde{G}_{L,2}^TD_1\widetilde{G}_{L,2} + (-\widetilde{G}_{L,1})^TD_2(-\widetilde{G}_{L,1})
\end{equation*} where 
 $\widetilde{G}_{L,l}= \left[ g_{e,v}^{L,l} \right]_{N_T \times \widetilde{N}_v} =\left[ \sqrt{\vert T \vert } \partial_l \varphi_v^L(x_T)  \right]_{N_T \times \widetilde{N}_v}$ and $\widetilde{C}_L=\left[ \begin{array}{c}
\widetilde{G}_{L,2} \\
-\widetilde{G}_{L,1}
\end{array} \right]_{2N_T\times \widetilde{N}_v}.$

Note that $(u_{CR},\widetilde{w}_L)$ satisfy the $2 \times 2$ uncoupled system 
\begin{equation*}
\left[ \begin{array}{cc}
A_{CR} & 0 \\ 0 & \widetilde{A}_L\\
\end{array} \right] \left[ \begin{array}{c}
u_{CR} \\ \widetilde{w}_L\\
\end{array} \right] = \left[ \begin{array}{c}
b_{CR} \\ \widetilde{b}_L
\end{array} \right].
\end{equation*}

Denote 
\begin{equation}\label{eq:A}
\widehat{A}=\left[ \begin{array}{cc}
A_{CR} & 0 \\ 0 & \widetilde A_L\\
\end{array} \right], ~~~\widehat{u}=\left[ \begin{array}{c}
u_{CR} \\ \widetilde{w}_L
\end{array} \right]
\mbox{ and } 
\widehat{b}=\left[ \begin{array}{c}
b_{CR} \\ \widetilde{b}_L
\end{array} \right]
\end{equation}
and introduce the matrices $H=[G_{CR} \ \ \widetilde{C}_L]$ and
\begin{equation}\label{eq:M}
M=H^TDH=\left[ \begin{array}{cc}
A_{CR} & G^T_{CR}D\widetilde{C}_L \\ \widetilde{C}_L^TDG_{CR} & \widetilde{A}_L
\end{array} \right].
\end{equation}
The preconditioned system is given by 
\begin{equation}\label{preconditioned}
M^{-1} \widehat{A} \widehat{u}= M^{-1}\widehat{b}.
\end{equation}
\begin{figure}\label{Mesh}
\centering
\includegraphics[scale=0.12]{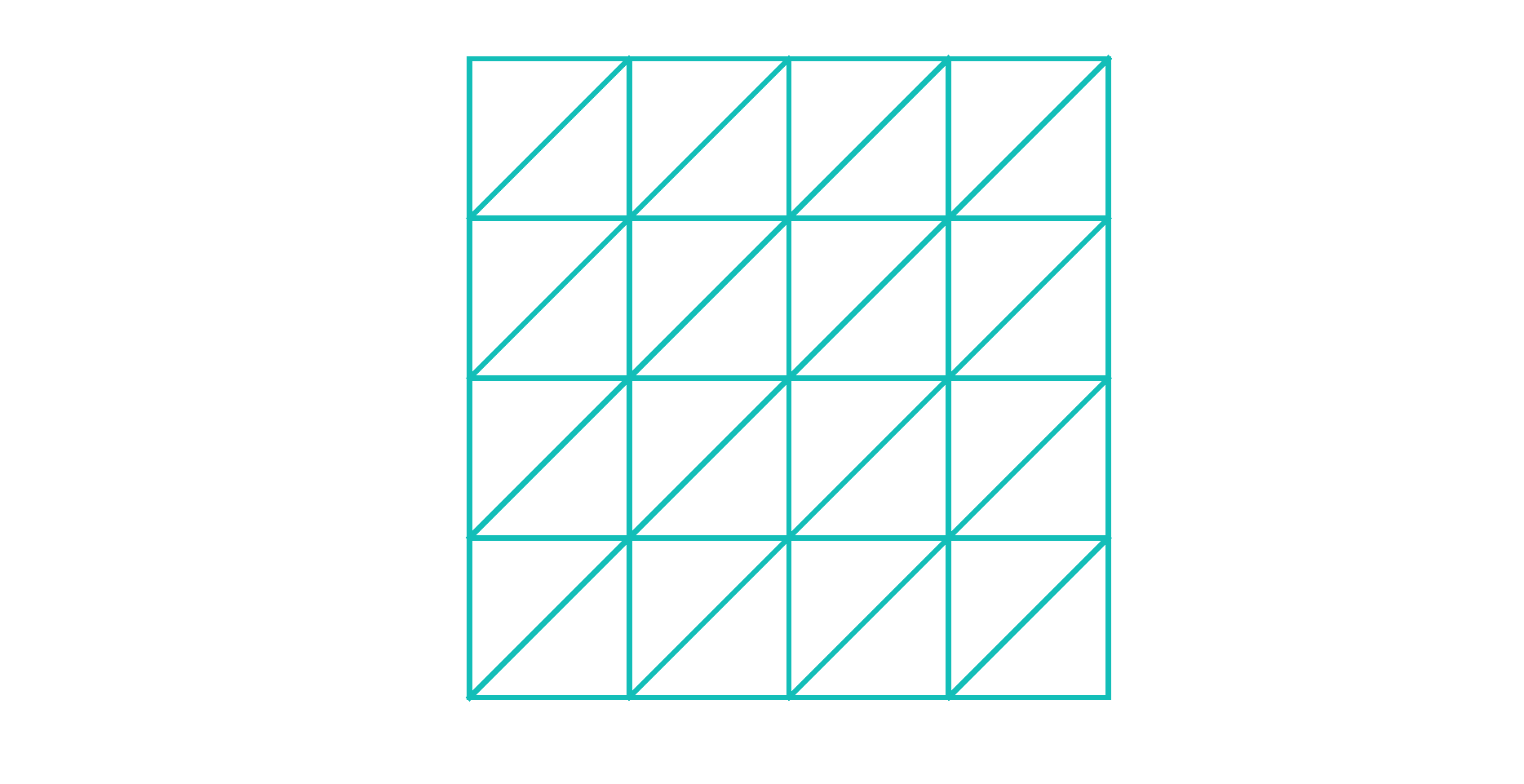}
\caption{Triangulation of $D=[0,1]^2$. }
\end{figure}

For any planar triangulation (with triangular elements) of a  simply connected domain  we have $2{N_T} = N_e +\widetilde{N}_v -1$ (where $N_e$ is the number of interior edges and $\widetilde{N}_v$
 is the number of vertices). See Figure \ref{Mesh} for the particular case of $\Omega=[0,1]^2$ and $\mathcal{T}^h$ constructed by dividing $\Omega$ into $n^2$ squares and further dividing each square into two triangles by adding and edge from the left-bottom vertex to right-top one. The following lemma shows that no extra computation is required to obtain basis of null spaces.
Also, recall that $G^TG$ is the stiffness 
matrix of the Laplace operator.\\

\noindent{\bf Lemma 1.} {\it We have
\begin{enumerate}
    \item[(a)] $H=[G_{CR} \ \ \widetilde{C}_L]$ is a square matrix of size $2N_T \times 2N_T$.
    \item[(b)] $G_{CR}^T  \widetilde{C}_L=0$. 
    \item[(c)] Because of (b), $H$ is non singular and $\widetilde{C}_L$ spans 
    the kernel of $G_{CR}^T$. Also $G_{CR}$ spans  the kernel of $\widetilde{C}_L^T$.
    \item[(d)] $M=H^TDH$ is the product of three square matrices. Therefore the solution of $M\widehat{v}=\widehat{r}$ can be computed as $\widehat{v}= H^{-1}D^{-1}H^{-T}\widehat{r}.$
\end{enumerate} }      
\noindent{\bf Proof}: We prove (b). Let $e$ be an interior edge and $v$ a vertex of $\mathcal{T}^h$. Then 
    \begin{eqnarray*}
    (G_{CR}^T \widetilde{C}_L)_{e,v} & = & (G_{CR,1}^T  \widetilde{G}_{L,2} - G_{CR,2}^T \widetilde{G}_{L,1})_{e,v}\\
    & = & \sum_{T \in \mathcal{T}^h} g_{e,T}^{CR,1} g_{v,T}^{L,2} - \sum_{T \in \mathcal{T}^h} g_{e,T}^{CR,2} g_{v,T}^{L,1}\\
    & = & \sum_{T \in \mathcal{T}^h} \vert T \vert \left[ \partial_1 \phi_e^{CR}(x_T) \partial_2 \varphi_v^L(x_T) - \partial_2 \phi_e^{CR}(x_T) \partial_1 \varphi_v^L(x_T) \right] \\ 
    & = & \sum_{T \in \mathcal{T}^h}  \vert T \vert \nabla \phi_e^{CR}(x_T) \cdot \overrightarrow{\mbox{curl }} \varphi_v^L(x_T) \\
    & = & \sum_{T \in \mathcal{T}^h} \int_T \nabla \phi_e^{CR}(x) \cdot \overrightarrow{\mbox{curl }} \varphi_v^L(x)\, dx \\
    & = & \sum_{T \in \mathcal{T}^h} \int_{\partial T}  \phi_e^{CR}(x) \,\overrightarrow{\mbox{curl }}  \varphi_v^L(x) \cdot \vec{\eta} \,  dx = 0.
    \end{eqnarray*}

We have the following condition number bound.  \\

\noindent{\bf Theorem 1.} {\it Let $\kappa_{\min} \leq  \kappa(x) \leq \kappa_{\max}$
    and denote $\eta = \kappa_{\max}/\kappa_{\min}$ the contrast. Then} 
    \[
\text{cond}\,(H^{-1}D^{-1}H^{-T} A) \leq 2 \eta -1.
    \]
        {\bf Proof}: let $s = u_{CR}^T A_{CR} u_{CR} + \widetilde{w}_{L}^T \widetilde{A}_{L} \widetilde{w}_{L}$, using Lemma 1 
      (b), the result follows from
$
  2|u_{CR}^T G_{CR}^T D \widetilde{C}_L \widetilde{w}_{L}| = 
  2|u_{CR}^T G_{CR}^T (D - D(k_{\min})) \widetilde{C}_L \widetilde{w}_{L}| \leq
  (1 - 1/\eta)s. 
$

\section{PCG 
for the block system}
We propose  to solve 
$
\widehat{A}\widehat{u}=\widehat{b}
$
with $\widehat{A}$ and $\widehat{b}$ defined in 
\eqref{eq:A} with $\widetilde{b}_L=0$ using  PCG 
with preconditioner $M$ in 
\eqref{eq:M}. See \eqref{preconditioned}. Recall that we use the construction in 
Section \ref{auxpro}.  For the  numerical test we compute an LU or QR factorizations for $H$ and apply $M^{-1}=H^{-1}D^{-1}H^{-T}$. Note that $M^{-1}$
depends on the coefficient $\kappa$ only through the diagonal matrix $D=D(\kappa)$ defined in \eqref{def:Dkappa}.
\subsection{Numerical tests for exponential covariance function}\label{application}
For problem (\ref{eq:hetdiff}) we consider the coefficient $\kappa$ of the form $\kappa(x,\omega)=e^{c(x,\omega)},$ where the stochastic process $c$ is defined by the Karhunen-Loève expansion with associated covariance function 
\begin{equation}\label{CE}
 \mathbf{c}( x,x' )=\exp \left( - \dfrac{1}{2} \parallel x-x' \parallel^2 \right).
\end{equation}
We approximate the expected value $\overline{u}(x)$ of the solution (\ref{eq:hetdiff}), through Monte Carlo method with $R$ realizations. In Table  1 we show the mean and variance of 
condition number of the preconditioned system, the number of iterations and the 
contrast $\max_{x}\kappa(x,\omega)/\min_{x} \kappa(x,\omega)$  during the Monte Carlo solve.
The small variance in the condition number indicates low dependence of the method on the parameter $\kappa$.

\begin{table}[H]\label{Tab2}
\begin{center}
\begin{tabular}{|l|r|r|r|}
\hline
 & \textbf{Condition} & \textbf{Iterations} & \textbf{Contrast}   \\
\hline 
\textbf{Mean} & 1.79 & 7.32 & 5.65 \\ \hline
\textbf{Variance} & 0.23 & 1.46 & 23.91 \\ \hline 
\end{tabular}
\caption{Condition number, number of iteration and coefficient contrast  in the CG method for the Monte Carlos computation of $\overline{u}(x)$ for (\ref{eq:hetdiff}). The log-coefficient $c$ is given by a truncated KL expansion with $K=15$ terms with covariance function shown in (\ref{CE}). We use $N=40$ elements in each direction and $R=1000$ realizations of the Monte Carlo method.}
\label{table_N40}
\end{center}
\end{table}

\subsection{Matérn class of covariance functions}

Now, the coefficient $\kappa$ is defined with the Matérn class of covariance functions 
\begin{equation}\label{Matern_Class}
    \mathbf{c}_{_{\mbox{Matern}}}(x,x')=\dfrac{2^{1-\nu}}{
    \Gamma(\nu)}\left(   \dfrac{\sqrt{2 \nu} \| x-x' \|}{l}\right)^\nu K_\nu \left(  \dfrac{\sqrt{2 \nu} \| x-x' \|}{l} \right)
\end{equation}
with (probabilistic) parameters $\nu,l>0$ and $K_{\nu}$ is the modified Bessel function of the second kind. With this function in the KL expansion, we obtain the results 
in Table 2 and 3. In Table 2 we verify converge of the Monte Carlo approximation and in Table 3 we show the dependence of the condition number, number of iteration and coefficient contrast. We note that the small variance of the number of iterations and  the value of the condition number indicate that the iteration do not depend much on the parameter $\kappa(x,\omega)$. Precise statements and results are object of current research and will be presented elsewhere.

\begin{table}[H]\label{Tab2}\small
\begin{center}
\begin{tabular}{|c|c|c|c|c|c|}
\hline
 R  & K=10 & K=20 & K=30 & K=40 & K=50  \\
\hline 
$10^3$ & 0.035066 & 0.034986 & 0.028177 & 0.028818 & 0.018414 \\ \hline
$10^4$ & 0.002926 & 0.000647 & 0.004533 & 0.004176 & 0.001530 \\\hline 
$10^5$ & 0.003007 & 0.001742 & 0.002240 & 0.002635 & 0.001134 \\\hline 
\end{tabular}
\caption{Error $H^1$ for the Monte Carlo 
approximation of $\overline{u}(x)$ where 
$u$ solves (\ref{eq:hetdiff}). The coefficient $c$ is a truncated KL expansion with $K$ terms constructed from covariance function shown in (\ref{Matern_Class}) with $\nu=0.5$ and $l=5$. We use $N=20$ elements in each direction and $R$ realization. The reference solution uses the same parameters  $\nu, l, N, R$ and $K=100$ terms of KL series.}
\vspace{-1cm}
\label{table_Matern_Class}
\end{center}
\end{table}
\begin{table}[H]\label{Tab2}
\begin{center}
\begin{tabular}{|c|c|c|c|}
\hline
 & \textbf{Condition} & \textbf{Iterations} & \textbf{Contrast}   \\
\hline 
\textbf{Mean} & 3.07 & 11.1 & 11.18 \\ \hline
\textbf{Variance} & 0.73 & 1.3 & 67.28 \\ \hline
\end{tabular}
\caption{Condition number, iterations numbers and contrast of 
coefficient $\kappa$  in the CG method in the Monte Carlo computation of $\overline{u}(x)$ solution of (\ref{eq:hetdiff}). The log-coefficient $c$ given as a truncated KL expansion with $K=30$ terms constructed from the covariance function shown in (\ref{Matern_Class}) with $\nu=0.5$ and $l=1$. We use $N=20$ elements in each direction and $R=1000$ realization of the Monte Carlo method.}
\label{table_MClassN20}
\end{center}
\end{table}
\vspace{-.5cm}
\noindent{\bf Acknowledgements.} The authors are grateful to Professor Zlatko Drmac from Univesity of Zagreb for introducing the third author to the natural
factor formulation of the stiffness matrices in finite element computations.

\bibliographystyle{abbrv}
\bibliography{references}

\end{document}